\documentclass[12pt]{article}
\usepackage{amsmath}
\usepackage{amssymb}
\begin{document}

\title{Finding and Investigating Exact Spherical Codes}
\author{Jeffrey Wang}

\begin{abstract}
In this paper we present the results of computer searches using a variation of an energy minimization algorithm used by Kottwitz for finding good spherical codes. We prove that exact codes exist by representing the inner products between the vectors as algebraic numbers. For selected interesting cases, we include detailed discussion of the configurations. Of particular interest are the $20$-point code in $\mathbb{R}^6$ and the $24$-point code in $\mathbb{R}^7$, which are both the union of two cross polytopes in parallel hyperplanes. Finally, we catalogue all of the codes we have found.
\end{abstract}

\maketitle

\section{Introduction}

Given $N$ points that lie on the unit sphere $S^{n-1}$ in $\mathbb{R}^n$, we wish to determine how they should be placed so that the minimal distance between any two points is maximized. Any set of points on the unit sphere is called a spherical code, and the problem of finding the best code has been proposed many times and in many contexts (though usually only in three dimensions), from packing circles on a sphere to distributing orifices on pollen-grains (e.g. \cite{clare}, \cite{tammes}). The only known optimal solutions in three dimensions are for $N \le 12$ and $N = 24$, while in four dimensions only the cases $N \le 8$, $N = 10$ \cite{410}, and $N = 120$ have been proven. The remaining known optimal codes are for $N \le 2n$ for any $n$, in which the points form either a simplex or a subset of the cross polytope; certain codes derived from the Leech lattice and the $E_8$ root system; and lastly an infinite family based on isotropic subspaces \cite{energy}. Excluding these few cases, the best known codes have been found either by specific constructions or, more commonly, by computer searches using various optimization algorithms. Given the difficulty of proving the optimality of even very small codes, most work related to this problem has been in finding close approximations to good configurations. The most extensive table of codes was created by N. J. A. Sloane, with the collaboration of R. H. Hardin, W. D. Smith and others, and is available electronically \cite{sloane}.

In this work we use the technique of energy minimization to find good spherical codes. Leech first observed the possibility of such an approach, and it has been used in previous work several times \cite{leech}. Kottwitz gave a fairly comprehensive list of three-dimensional codes for up to $90$ points \cite{kottwitz}, which was expanded later by Buddenhagen and Kottwitz while searching for codes with multiplicity greater than one (i.e. codes for which there exist distinct configurations that obtain the same optimal distance) \cite{budden}. Nurmela investigated some global optimization methods based on energy minimization and provided numerical results for codes in up to five dimensions \cite{nurmela}, but the most successful implementations have been based on using a large number of random starts and a local optimization algorithm so that there is a high probability that at least one of them converges to the global optimum. Our algorithm has several changes to improve on old implementations as well, such as reducing exponent bias and choosing a different local optimization method. This has resulted in three improved spherical codes in four dimensions as well as new higher dimensional codes, particularly in six and seven dimensions, that exhibit interesting symmetries. Section \ref{sec:method} describes our algorithm in detail and Section \ref{sec:improved} gives a brief analysis of the improved codes.

The remainder of this work focuses on providing rigorous analysis of codes, particularly showing that the best known codes can be represented as exact configurations in terms of algebraic numbers as opposed to close estimates. This is an important step toward rigorizing computer solutions to the problem and allows us to observe true equalities and relations between points and edges in the code. Buddenhagen and Kottwitz did similar work while looking for three-dimensional codes with multiplicity greater than one, and provided detailed exact descriptions of the two distinct optimal $15$-point codes \cite{budden}. In several cases we were able to identify a considerable amount of underlying symmetry and structure using the methods of Section \ref{sec:identify}. For three particularly interesting cases we provide brief discussion of the configurations in Section \ref{sec:selected}. The final section gives tables of the exact codes, in dimensions four through eight, based on their minimal polynomials.

\section{Methodology}
\label{sec:method}
We consider an inverse power law potential function on the spherical code $C = \{x_1, x_2, \ldots, x_N\} \subset S^{n-1}$ defined by
\begin{equation}
\displaystyle E = \sum_{1 \le i<j \le N} \left(\frac{\alpha}{|x_i-x_j|}\right)^s,
\end{equation}
where $\alpha$ is a constant to prevent overflow, $s$ is the exponent of the inverse power law, and $|x_i-x_j|$ is the Euclidean distance between the two points in $\mathbb{R}^n$. As $s \rightarrow \infty$, the smallest distance will be the dominating term in the energy expression. Hence, minimizing this potential function over all possible $C$ will give approximate codes that get better as $s$ increases. In certain highly symmetrical cases, when $s$ is sufficiently large, the method will produce the actual optimal code, whereas in general it converges to the optimal code as $s \rightarrow \infty$.

Several issues arise, however, if one tries to minimize $E$ using a large $s$ with a random configuration of points. There is an implementation problem, in that $\alpha$ cannot be chosen well so that $E$ neither overflows nor is too small. It has no abstract mathematical role, but it is important to consider when dealing with floating-point arithmetic. If a poor choice of $\alpha$ is used with a power of $s = 1,000,000$ for example, a slight deviation in $|x_i-x_j|$ from $\alpha$ by only $0.001$ will result in either $1.001^{1000000} \approx 10^{434}$ or $0.999^{1000000} \approx 10^{-435}$, neither of which would be handled well by a computer. There is also a problem regarding the behavior of the function $E$ as $s$ increases because the number of local minima increases as well. Starting with a high exponent in the inverse power interaction will usually converge quickly to a poor local minimum.

Taking into consideration these problems, Kottwitz used the following algorithm \cite{kottwitz}:
\newcounter{count1}
\begin{list}{\arabic{count1}.}{\usecounter{count1} \setlength{\rightmargin}{\leftmargin}}
\item Start with $s = 80$ and a random configuration of points;
\item Run a local optimization algorithm until convergence;
\item Double $s$ and repeat the previous step starting from the local optimum found, stopping after the optimization is run for $s = 1,310,720$.
\end{list}
This remedies the issues described because $\alpha$ can be chosen before each optimization to reflect the minimal distance of the configuration, and the points will tend toward a good code without getting stuck early on in a bad local minimum. For our work we use a similar strategy, but to avoid the bias introduced by fixing an initial exponent we start with a random exponent between $10$ and $160$ each time. Furthermore, we allow the doubling procedure to continue until the exponent exceeds $100,000,000$. Before each optimization procedure, we choose $\alpha$ to be the exact minimal distance between all pairs of points, as done by Nurmela \cite{nurmela}. This is a good choice because as the power increases the improvement in the minimal distance decreases, and we have found that this balance keeps the energy expression in the desired range.

While Kottwitz used gradient descent and Newton-Raphson methods to optimize the energy function at each step \cite{budden}, \cite{kottwitz}, we picked the nonlinear conjugate gradient algorithm for several reasons. Firstly, it is easy to implement yet much more effective than gradient descent. It also only requires the gradient at each step and does not need to compute the Hessian, which is costly and very complex especially as the dimension increases. We experimented with three versions of the nonlinear conjugate gradient algorithm: Fletcher-Reeves, Polak-Ribi\`{e}re, and Hestenes-Stiefel. We found that the Hestenes-Stiefel formula for updating the conjugate direction was most effective for our model, converging faster than Polak-Ribi\`{e}re and displaying better results than Fletcher-Reeves.

We must also consider the constrained nature of optimization on a sphere, as the conjugate gradient algorithm does not accommodate for constraints on the domain. Our solution to this is to, after each movement in the conjugate direction, simply scale each point in $C$ back onto the sphere. We realize, however, that as the code approaches a local minimum, the points will tend to move almost directly away from the sphere because most of the movement along the tangent space will be balanced by their neighbors. When scaled they will consequently appear to have moved very little if at all, slowing convergence and possibly even getting stuck before converging. To compensate for this, when the exponent in the inverse power law interaction is high, for each point we only look at the component of the gradient that lies in the tangent space to the sphere at that point. We do not implement such a procedure for lower exponents because the conjugate directions depend on previous gradients as well so the magnitudes should be consistent relative to each other.

Using the method that has been described, implemented in C++, we ran the optimization with at least $1000$ random starts. Taking the best result of these runs, we then applied a similar optimization to the configuration to improve the precision using a starting exponent of $10,000,000$ and doubling until $640,000,000$. Once this converged we were able to obtain fairly accurate results for our codes, but not accurate enough to find their exact representations. Finally, we found the set of shortest edges to approximately five or six decimal places and set them equal to each other, which determined a rigid structure solved using Newton-Raphson methods for a small number of iterations to a precision of about $500$ decimal digits. Then we placed the remaining points, the rattlers, by moving them iteratively as far as possible from the rigid points, again maximizing the minimal distance. We used the PARI/GP software for the last two optimization steps.

After the full procedure, our set of codes performed very well against previous data. For the three-dimensional codes we ran our program on, we reproduced the best known results in every case, while for four-dimensional codes we were able to find three codes with $40$, $68$, and $71$ points that improved on the best previous results. In higher dimensions for smaller numbers of points we are confident that we were able to find good (and most likely optimal) configurations.

\section{Improved Codes in $\mathbb{R}^4$}
\label{sec:improved}
Our work has improved on three codes in four dimensions, namely those with $40$, $68$, and $71$ points. They are briefly examined in this section. Note that the improvements are small -- on the order of $1/1000$th or less. A catalogue of all of the codes can be found online at
\begin{center}
\texttt{http://www.aimath.org/data/paper/WangSphericalCode/}.
\end{center}
The data available online reflects only results prior to the optimization done in PARI/GP.

\subsection{$40$ points}

The maximal cosine of $0.65049780106271773133\ldots$ occurs $108$ times, but is the only inner product that occurs more than twice, suggesting that this configuration is quite asymmetric.

\subsection{$68$ points}

The maximal cosine of $0.75104449257228207352\ldots$ occurs $196$ times, and there is one other inner product - $0.74925795575260304153\ldots$ - that occurs three times. Every other edge length appears to be unique.

\subsection{$71$ points}

The maximal cosine of $0.75637601134814761871\ldots$ occurs $199$ times, and there are three other inner products that occur three times each in the code; they are $0.75152318440020477595\ldots$, $0.75574916415250115097\ldots$, and $0.75632115855377282465\ldots$. The other edges are unique.

\section{Properties for Identifying Codes}
\label{sec:identify}
We examine some properties of spherical codes that can hopefully characterize them more substantially and rigorously than a list of floating-point coordinates. First, we introduce the $N \times N$ Gram matrix $G$ such that $G_{ij} = \langle x_i, x_j \rangle$, where $\langle x, y \rangle$ denotes the standard inner product between the unit vectors $x$ and $y$. Since all vectors on the unit sphere have magnitude one, we also have $G_{ij} = \cos{(\theta_{ij})}$, where $\theta_{ij}$ is the angle between the vectors $x_i$ and $x_j$.

\subsection{Edge Length Occurrences}

One effective way of determining whether or not a spherical code has symmetry is to look at the number of occurrences of each edge length or the number of unique edge lengths between distinct points, which we will call $m$. Since the distance and inner product are related by
\begin{equation}
|x_i-x_j|^2 = 2-2\langle x_i, x_j \rangle,
\end{equation}
this consists of simply looking at the unique entries of the Gram matrix and the number of times each is repeated. We note that considering edge length occurrences is a generic and flexible strategy for describing codes and detecting overall symmetry, especially when $m$ is small relatively to the code size. To find more specific symmetries, we can look at the edge lengths emanating from any one point in the code and check if any other points shared these. In particular, we can often group vertices together as equivalent if the lists of edge lengths emanating from them are exactly the same. We note that this is a necessary but not sufficient condition for the vertices to be equivalent under the action of the symmetry group of the configuration.

\subsection{Exact Algebraic Numbers}

One inevitable concern that arises from computer searches for spherical codes is whether or not the set of points found by the computer really exists. Numerical results can be useful for recognizing the exact solutions, but with finite precision it is not clear if all the distances are really what we think they are. Consequently it is important to have some sort of method of converting the work of computer approximations into definitive arrangements so that they may be analyzed with more certainty and have hope of being proven optimal as well.

Since a code is the solution to a number of polynomial equalities between the shortest edges, the coordinates of each rigid point in the code are algebraic and lie in a number field. Moreover, we then know that the inner products are algebraic as well. Once we have determined the coordinates of the points in the code to high precision we are able to effectively compute, using the PARI/GP math software command \texttt{algdep}, the minimal polynomial for the maximal inner product between two distinct vectors, i.e.,
\begin{equation}
\displaystyle u = \max_{1 \le i < j \le N} \langle x_i, x_j \rangle.
\end{equation}
This was done for a number of three-dimensional codes by Buddenhagen and Kottwitz, and we extend their work. Just looking for a minimal polynomial $P$ so that $P(u) \approx 0$ can be misleading, however, because even with high precision it could be the case that $u$ is in fact not the root of $P$, but just very close. On the other hand, if we could show that the code with exactly the root of $P$ (closest to $u$, to be specific) as the maximal inner product does exist, we would have an exact arrangement that is practically equivalent. We have been able to do so for most codes that we could calculate the minimal polynomial of.

To verify the existence of the codes, we first use the PARI/GP command \texttt{lindep} to specify the elements of $G$ in terms of polynomials in some primitive element of the number field. Empirically we have found that $u$ is often sufficient for doing this, though in several cases this has not been true (in particular for $22$ points in three dimensions, $21$ points in five dimensions, and $14$ points in six dimensions). Once we have the exact Gram matrix, we may check if $G$ is indeed a Gram matrix of a set of points on the unit sphere in $\mathbb{R}^n$, namely that it satisfies the following properties: (a) $G$ is symmetric; (b) the diagonal elements of $G$ are all equal to $1$; (c) $G$ has rank at most $n$; and (d) $G$ is positive semidefinite. If $G$ satisfies each of these conditions, it must be the Gram matrix of a spherical code, and thus we are able to prove the existence of such a code. The first two properties are easy to check, and the last two can be tested by observing the characteristic polynomial $p_G(\lambda)$ of $G$, which will be of degree $N$. Let
\begin{equation}
p_G(\lambda) = a_N\lambda^N+a_{N-1}\lambda^{N-1}+\cdots+a_1\lambda+a_0.
\end{equation}
It is known that $G$ has rank at most $n$ if and only if $p_G(\lambda)$ has a root at $\lambda = 0$ with multiplicity at least $N-n$, or equivalently $a_0 = a_1 = \cdots = a_{N-n-1} = 0$. Also, $G$ is positive semidefinite if and only if the remaining coefficients are nonzero and alternate in sign. This makes it relatively easy to check if $G$ is indeed a Gram matrix.

\section{Discussion of Selected Cases}
\label{sec:selected}
Here we present more detailed results of three cases which we found particularly interesting and which have not been previously described in the literature (except the last). Some brief discussion of other cases can be found online at \texttt{http://www.aimath.org/data/paper/WangSphericalCode/}.

\subsection{$12$ points in $\mathbb{R}^5$}

This is a very simple configuration, so it is not surprising that the Gram matrix is quite simple as well, consisting of only five unique entries (excluding the $1$'s along the diagonal). It is shown here:

\[ \left( \begin{array}{cccccccccc|cc}
1 & u & V & V & u & u & u & u & u & W & u & -u \\
u & 1 & u & V & V & u & u & u & W & u & u & -u \\
V & u & 1 & u & V & u & u & W & u & u & u & -u \\
V & V & u & 1 & u & u & W & u & u & u & u & -u \\
u & V & V & u & 1 & W & u & u & u & u & u & -u \\
u & u & u & u & W & 1 & V & u & u & V & -u & u \\
u & u & u & W & u & V & 1 & V & u & u & -u & u \\
u & u & W & u & u & u & V & 1 & V & u & -u & u \\
u & W & u & u & u & u & u & V & 1 & V & -u & u \\
W & u & u & u & u & V & u & u & V & 1 & -u & u \\
\hline
u & u & u & u & u & -u & -u & -u & -u & -u & 1 & -1 \\
-u & -u & -u & -u & -u & u & u & u & u & u & -1 & 1 \\
\end{array} \right), \] \\
where $V$ and $W$ satisfy $V > W$ and are given in terms of $u$ by
\begin{center}
$\begin{array}{ccccccc}
V &=& \frac{5}{2}u^2-u-\frac{1}{2}, & & W &=& -5u^2-4u. \\
&
\end{array}$
\end{center}

There are two aspects of this Gram matrix that are particularly intriguing. The first one that stands out is that all the edges between $V_1$ and $V_2$ correspond to an inner product of $u$ or $-u$, which is a result of the points in $V_2$ being antipodal. The second is that every element is represented as a polynomial in $u$ of degree two or less. The minimal polynomial, however, is $25u^4+30u^3+24u^2+2u-1$, which has degree four. This is unusual because in most cases a significant number of elements of the Gram matrix are polynomials in $u$ of maximal degree, i.e. one less than the degree of the minimal polynomial. While we do not have an explanation for this, it may be significant to the structure of the code.

\subsection{$20$ points in $\mathbb{R}^6$}

This configuration, the first especially nice one, is highly symmetric; each point in the entire configuration is equivalent to every other one. Each point is at the minimal distance with $11$ other points in the code. Furthermore, since the minimal polynomial is $14u-3$, all of the entries in the Gram matrix are rational numbers, specifically multiples of $\frac{1}{14}$. The dot products corresponding to the edges that emanate from any point are (in order from largest to smallest)
\begin{equation}
\left\{1, \frac{3}{14}, 0, -\frac{1}{14}, -\frac{5}{14}, -\frac{3}{7}, -\frac{4}{7}, -\frac{9}{14}\right\}.
\end{equation}
The only repetitions are $\displaystyle \frac{3}{14}$, which shows up $11$ times, and $\displaystyle -\frac{9}{14}$, which shows up three times.

It turns out that this code is the union of two five-dimensional cross polytopes in parallel hyperplanes. We can orient them with respect to one another: if we choose coordinates to allow one of the cross-polytopes to be given by plus or minus the standard orthonormal basis vectors, then after projection into the equatorial hyperplane and rescaling, the other cross polytope is given by

\begin{eqnarray*}
& \displaystyle \pm \left(\frac{2}{11}, -\frac{6}{11}, -\frac{6}{11}, \frac{3}{11}, -\frac{6}{11}\right) \\
& \displaystyle \pm \left(\frac{6}{11}, \frac{2}{11}, \frac{6}{11}, \frac{6}{11}, -\frac{3}{11}\right) \\
& \displaystyle \pm \left(-\frac{3}{11}, -\frac{6}{11}, \frac{2}{11}, \frac{6}{11}, \frac{6}{11}\right) \\
& \displaystyle \pm \left(\frac{6}{11}, \frac{3}{11}, -\frac{6}{11}, \frac{2}{11}, \frac{6}{11}\right) \\
& \displaystyle \pm \left(\frac{6}{11}, -\frac{6}{11}, \frac{3}{11}, -\frac{6}{11}, \frac{2}{11}\right) \\
\end{eqnarray*}

The way the orientation is determined is still unknown, but there is evidence of some sort of underlying structure that is worth investigating further.

\subsection{$24$ points in $\mathbb{R}^7$}

This is the other particularly nice configuration. All $24$ points in this code are symmetric with respect to one another and the minimal polynomial is $19u^2+2u-1$, so the Gram matrix has very few entries and they are all simple. The only entries other than the $1$'s along the diagonal are $u$, $-u$, $-3u$, and $2u-1$. Each vector has inner product $u$ with $15$ other vectors, $-u$ with $2$ other vectors, $-3u$ with $5$ other vectors, and $2u-1$ with a single other vector. This configuration appears to be a very good code as both the $22$- and $23$-point codes are the same as this with points removed.

It is very similar to the $20$-point code in $\mathbb{R}^6$ as it also consists of two cross-polytopes in parallel hyperplanes. As with the other case, we can orient the cross polytopes relative to each other by setting one to be plus or minus the standard orthonormal basis vectors. After projection into the equatorial hyperplane and rescaling, the other cross polytope is then given by the following vectors:

\begin{eqnarray*}
& \displaystyle \pm \left(0, \frac{1}{\sqrt{5}}, -\frac{1}{\sqrt{5}}, -\frac{1}{\sqrt{5}}, \frac{1}{\sqrt{5}}, \frac{1}{\sqrt{5}}\right) \\
& \displaystyle \pm \left(\frac{1}{\sqrt{5}}, 0, \frac{1}{\sqrt{5}}, -\frac{1}{\sqrt{5}}, -\frac{1}{\sqrt{5}}, \frac{1}{\sqrt{5}}\right) \\
& \displaystyle \pm \left(-\frac{1}{\sqrt{5}}, \frac{1}{\sqrt{5}}, 0, \frac{1}{\sqrt{5}}, -\frac{1}{\sqrt{5}}, \frac{1}{\sqrt{5}}\right) \\
& \displaystyle \pm \left(-\frac{1}{\sqrt{5}}, -\frac{1}{\sqrt{5}}, \frac{1}{\sqrt{5}}, 0, \frac{1}{\sqrt{5}}, \frac{1}{\sqrt{5}}\right) \\
& \displaystyle \pm \left(\frac{1}{\sqrt{5}}, -\frac{1}{\sqrt{5}}, -\frac{1}{\sqrt{5}}, \frac{1}{\sqrt{5}}, 0, \frac{1}{\sqrt{5}}\right) \\
& \displaystyle \pm \left(\frac{1}{\sqrt{5}}, \frac{1}{\sqrt{5}}, \frac{1}{\sqrt{5}}, \frac{1}{\sqrt{5}}, \frac{1}{\sqrt{5}}, 0\right) \\
\end{eqnarray*}

This is a very structured set of vectors, but also differs from the vectors in the six dimensional case. It would be an interesting question to determine how in general the cross polytopes should be arranged relative to each other to maximize the minimal distance.

\section{Catalogue of Exact Configurations}
\label{sec:catalogue}
We present several tables of results for newly calculated exact spherical codes in dimensions four through eight, specifically the maximal inner product $u$ and the minimal polynomial. In each dimension, we attempted to calculate a minimal polynomial for the maximum inner product of the configurations for up to $24$ points ($27$ in four and six dimensions). We were able to do so for quite a few, but there are a still a number of codes that we were unable to do this for. In general we attempted this calculation using $500$ decimal places of the maximum inner product.

Moreover, for the minimal polynomials we calculated, we were able to verify all but two of them using the methods of Section \ref{sec:identify}. These are $21$ points in five dimensions and $14$ points in six dimensions, both of which have polynomials of high degree. The case of $27$ points in four dimensions should also be noted because it has three rattlers, and since their positions are not fixed relative to the other points we did not include them in our verification of the minimal polynomial.

Lastly, we have a single new minimal polynomial in three dimensions, concerning the $22$-point code. The minimal polynomial we found is $486u^{18} + 13113u^{17} + 114996u^{16} + 117476u^{15} + 658256u^{14} + 378752u^{13} - 347056u^{12} - 121388u^{11} + 81724u^{10} - 70886u^9 - 55992u^8 + 12716u^7 + 6528u^6 - 2392u^5 - 208u^4 + 284u^3 + 14u^2 + 5u + 4$. But as in the two cases mentioned previously, we have been unable to verify this through finding an exact Gram matrix.

Note that Tables \ref{tab:4d} and \ref{tab:5d} are the same codes as recorded in \cite{sloane} but with minimal polynomials associated with them. Also, some codes predate references to \cite{ericson}, but it is a convenient, systematic reference work.

\begin{table}[ht]
\caption{Four-Dimensional Codes}
{\footnotesize
\begin{tabular}{|c|c|c|c|}
\hline
$N$ & $u$ ($20$ decimal places) & Minimal Polynomial & Ref \\
\hline
$9$ & $0.16201519961163454918$ & $16u^3-16u^2-4u+1$ & \cite{ericson} \\
\hline
$10$ & $0.16666666666666666667$ & $6u-1$ & \cite{ericson} \\
\hline
$11$ & $0.23040556359455544174$ & $8u^3-12u^2-2u+1$ & \\
\hline
$12$ & $0.25000000000000000000$ & $4u-1$ & \cite{ericson} \\
\hline
$13$ & $0.30729565398102882233$ & $5632u^9+9472u^8-3072u^7-5888u^6+544u^5+$ & \\
& & $944u^4+152u^3-44u^2-14u-1$ & \\
\hline
$14$ & $0.31951859421260363550$ & $58u^7+174u^6+140u^5-16u^4-$ & \\
& & $54u^3-10u^2+4u+1$ & \\
\hline
$15$ & $0.35099217594534630330$ & $36u^4-18u^3+10u^2-1$ & \\
\hline
$16$ & $0.38762817712253427776$ & $256u^{10}+1024u^9+256u^8-1152u^7+160u^6+$ & \\
& & $176u^5-144u^4+36u^3+37u^2-6u-3$ & \\
\hline
$17$ & $0.41226632322755925382$ & Unknown & \\
\hline
$18$ & $0.42281941407305934403$ & $14424u^{16}+42932u^{15}+18232u^{14}-62100u^{13}-$ & \\
& & $53831u^{12}+41528u^{11}+46442u^{10}-18248u^9-$ & \\
& & $20977u^8+6180u^7+5372u^6-1556u^5-$ & \\
& & $721u^4+240u^3+34u^2-16u+1$ & \\
\hline
$19$ & $0.43425854591066488219$ & $3u^2+u-1$ & \\
\hline
$20$ & $0.43425854591066488219$ & $3u^2+u-1$ & \cite{ericson} \\
\hline
$21$ & $0.47138085850731791682$ & $16u^8-128u^7-64u^6+32u^5+$ & \\
& & $72u^4+32u^3-16u^2-8u+1$ & \\
\hline
$22$ & $0.49788413084355235629$ & Unknown & \\
\hline
$23$ & $0.50000000000000000000$ & $2u-1$ & \\
\hline
$24$ & $0.50000000000000000000$ & $2u-1$ & \cite{ericson} \\
\hline
$25$ & $0.53731605665507787660$ & Unknown & \\
\hline
$26$ & $0.54078961769753707673$ & $3392u^6+2112u^5-496u^4-$ & \\
& & $656u^3-132u^2+6u-1$ & \\
\hline
$27$ & $0.55759135118017018253$ & $794u^5+393u^4-344u^3-82u^2+6u+1$ & \\
\hline
\end{tabular}
}
\label{tab:4d}
\end{table}

\begin{table}[ht]
\caption{Five-Dimensional Codes}
\centering
{\footnotesize
\begin{tabular}{|c|c|c|c|}
\hline
$N$ & $u$ ($20$ decimal places) & Minimal Polynomial & Ref \\
\hline
$11$ & $0.13285354259858991809$ & $45u^3-25u^2-5u+1$ & \cite{ericson} \\
\hline
$12$ & $0.15393160503302123095$ & $25u^4+30u^3+24u^2+2u-1$ & \\
\hline
$13$ & $0.18725985188285358702$ & $17u^3-5u^2-5u+1$ & \\
\hline
$14$ & $0.20000000000000000000$ & $5u-1$ & \\
\hline
$15$ & $0.20000000000000000000$ & $5u-1$ & \\
\hline
$16$ & $0.20000000000000000000$ & $5u-1$ & \cite{ericson} \\
\hline
$17$ & $0.27047583526857362209$ & $9u^4-16u^3-10u^2+1$ & \\
\hline
$18$ & $0.27550174165981923839$ & $484u^5-488u^4+97u^3+17u^2-u-1$ & \\
\hline
$19$ & $0.29182239902449014615$ & $57u^6-38u^5-109u^4+32u^3+23u^2-10u+1$ & \\
\hline
$20$ & $0.29938569289912478230$ & $5u^3+13u^2-u-1$ & \\
\hline
$21$ & $0.31491695717530346285$ & $869312u^{14}+8798656u^{13}-1062776u^{12}-$ & \\
& & $10586775u^{11}-968269u^{10}+3532907u^9+$ & \\
& & $188249u^8-659974u^7-11746u^6+72246u^5-$ & \\
& & $806u^4-5267u^3-97u^2+207u+21$ & \\
\hline
$22$ & $0.35499503416625620683$ & Unknown & \\
\hline
$23$ & $0.36977269694307633377$ & Unknown & \\
\hline
$24$ & $0.37423298246516725173$ & $1620u^{10}+5508u^9-5751u^8-2406u^7+2055u^6+$ & \\
& & $276u^5+559u^4-210u^3-127u^2+48u-4$ & \\
\hline
\end{tabular}
}
\label{tab:5d}
\end{table}

\begin{table}[ht]
\caption{Six-Dimensional Codes}
\centering
{\footnotesize
\begin{tabular}{|c|c|c|c|}
\hline
$N$ & $u$ ($20$ decimal places) & Minimal Polynomial & Ref \\
\hline
$13$ & $0.11307975214744721385$ & $96u^3-36u^2-6u+1$ & \cite{ericson} \\
\hline
$14$ & $0.13249092032347031437$ & $237299200u^{13}+463738880u^{12}+366062080u^{11}+$ & \\
& & $142462784u^{10}+23021376u^9-2398832u^8-$ & \\
& & $1849600u^7-343428u^6-1802u^5+11443u^4+$ & \\
& & $1390u^3-128u^2-20u+1$ & \\
\hline
$15$ & $0.14494897427831780982$ & $20u^2+4u-1$ & \\
\hline
$16$ & $0.17114764942939365334$ & $4000000u^8+15120000u^7+9896400u^6-$ & \\
& & $2451600u^5-4503600u^4-1511460u^3-$ & \\
& & $87480u^2+43740u+6561$ & \\
\hline
$17$ & $0.18327433702314481858$ & $400u^4+240u^3+16u^2-8u-1$ & \\
\hline
$18$ & $0.19781218860197545241$ & $784u^8-280u^7-1191u^6+120u^5+$ & \\
& & $218u^4-62u^3-31u^2+2u+1$ & \\
\hline
$19$ & $0.20022602589120548304$ & Unknown & \\
\hline
$20$ & $0.21428571428571428571$ & $14u-3$ & \\
\hline
$21$ & $0.24285284801369170251$ & Unknown & \\
\hline
$22$ & $0.24886966078882474478$ & Unknown & \\
\hline
$23$ & $0.25000000000000000000$ & $4u-1$ & \\
\hline
$24$ & $0.25000000000000000000$ & $4u-1$ & \\
\hline
$25$ & $0.25000000000000000000$ & $4u-1$ & \\
\hline
$26$ & $0.25000000000000000000$ & $4u-1$ & \\
\hline
$27$ & $0.25000000000000000000$ & $4u-1$ & \cite{ericson} \\
\hline
\end{tabular}
}
\label{tab:6d}
\end{table}

\begin{table}[ht]
\caption{Seven-Dimensional Codes}
\centering
{\footnotesize
\begin{tabular}{|c|c|c|c|}
\hline
$N$ & $u$ ($20$ decimal places) & Minimal Polynomial & Ref \\
\hline
$15$ & $0.09870177627236447933$ & $175u^3-49u^2-7u+1$ & \cite{ericson} \\
\hline
$16$ & $0.11332087960014474125$ & $85823999u^{12}+153503984u^{11}+163426326u^{10}+$ & \\
& & $92426704u^9+29708081u^8+9292960u^7+$ & \\
& & $2252404u^6+410976u^5+29489u^4-$ & \\
& & $3216u^3-746u^2-48u-1$ & \\
\hline
$17$ & $0.12484758986639552862$ & Unknown & \\
\hline
$18$ & $0.12613198362288317392$ & $47u^2+2u-1$ & \\
\hline
$19$ & $0.15659738541709551030$ & $1280u^6-352u^4-48u^3+37u^2+3u-1$ & \\
\hline
$20$ & $0.16952084719853722593$ & $23u^2+2u-1$ & \cite{ericson} \\
\hline
$21$ & $0.18152396080041583541$ & Unknown & \\
\hline
$22$ & $0.18274399763155681015$ & $19u^2+2u-1$ & \\
\hline
$23$ & $0.18274399763155681015$ & $19u^2+2u-1$ & \\
\hline
$24$ & $0.18274399763155681015$ & $19u^2+2u-1$ & \\
\hline
\end{tabular}
}
\label{tab:7d}
\end{table}

\begin{table}[ht]
\caption{Eight-Dimensional Codes}
\centering
{\footnotesize
\begin{tabular}{|c|c|c|c|}
\hline
$N$ & $u$ ($20$ decimal places) & Minimal Polynomial & Ref \\
\hline
$17$ & $0.08773346332333854567$ & $288u^3-64u^2-8u+1$ & \cite{ericson} \\
\hline
$18$ & $0.09946957270878709386$ & Unknown & \\
\hline
$19$ & $0.11140997502603998543$ & Unknown & \\
\hline
$20$ & $0.11949686668719356518$ & $244u^3+60u^2-2u-1$ & \\
\hline
$21$ & $0.13060193748187072126$ & $28u^2+4u-1$ & \\
\hline
$22$ & $0.13060193748187072126$ & $28u^2+4u-1$ & \cite{ericson} \\
\hline
$23$ & $0.15745541772761612286$ & Unknown & \\
\hline
$24$ & $0.15769214493936087799$ & $4u^4-4u^3-27u^2-2u+1$ & \\
\hline
\end{tabular}
}
\label{tab:8d}
\end{table}

\section{Acknowledgments}

I would like to thank my mentor Henry Cohn for suggesting this project and constantly providing invaluable assistance along the way. I would also like to thank the Microsoft Research Internship and Microsoft High School Internship programs for arranging my internship and allowing this all to happen. Lastly, I would like to thank the American Institute of Mathematics for hosting my data online.

\end{document}